\documentstyle[11pt]{article}

\newtheorem{th}{Theorem}[section]

\newtheorem{Proposition}[th]{Proposition}

\newtheorem{Corollary}[th]{Corollary}

\newtheorem{ccote}[th]{}

\newcommand{\cqfd}{\unskip\kern 6pt\penalty 500
\raise -2pt\hbox{\vrule\vbox to10pt{\hrule width
 4pt\vfill\hrule}\vrule}\smallskip}

\newcommand{\bbr}{{\bf R}}
\newcommand{\bbc}{{\bf C}}
\newcommand{\bbz}{{\bf Z}}
\newcommand{\bbq}{{\bf Q}}

\newcommand{\pcirc}{\kern .7pt {\scriptstyle \circ} \kern 1pt}
\newcommand{\iso}{\cong}

\newcommand{\hfl}[1]{\buildrel{#1}\over{\longrightarrow}}
\newcommand{\eqref}[1]{(\ref{#1})}

\title{Cohomology rings of symplectic cuts}
\author{Jean-Claude HAUSMANN
\and Allen KNUTSON
\footnote{Both authors thank the Fonds National Suisse de la
Recherche Scientifique for its support.}  }
\date{October 8, 1998}

\begin{document}    \maketitle
\begin{abstract}
The ``symplectic cut''  construction \cite{Le} produces two symplectic
orbifolds $C_-$ and $C_+$ from a symplectic manifold $M$ with a Hamiltonian
circle action. We compute the rational cohomology ring of $C_+$ in terms of
those of $M$ and $C_-$.
\end{abstract}

\section{Statement of the results} \label{intro}

Let $M$ be a symplectic $n$-manifold endowed with an Hamiltonian
$S^1$-action with moment map $f: M\to\bbr$. Suppose that $0$
is a regular value of $f$ and
consider the manifolds-with-boundary
$$M_-:=f^{-1}(\bbr_{\leq 0}) \quad , \quad
M_+:=f^{-1}(\bbr_{\geq 0}) \quad , \quad
M_0:=f^{-1}(0) = M_-\cap M_+ .$$
The symplectic cutting of $(M,f)$ at $0$, introduced by E. Lerman [Le],
can be  described as follows.
The $S^1$-action restricted to $M_0$ gives rise to the equivalence
relation on $M$:
$$x\sim y \ \Longleftrightarrow \
\left\{\begin{array}{lcl}
x\in M_0 & \hbox{ and } &
y=\theta\cdot x \hbox{ for some } \theta\in S^1 \\ \qquad \hbox{ or } \\
x=y & \hbox{ if }
x\notin
M_0.\end{array}\right.$$
As $0$ is a regular value of $f$, the stabilizers of points in $M_0$ are
finite cyclic subgroups of $S^1$.
Therefore, the following quotients are orbifolds:
$$%\begin{array}{lcl}
C_- :=  M_-\big/\sim  \quad , \quad
C_+ :=  M_+\big/\sim \quad , \quad
C_0  :=  M_0\big / S^1 %= C_-\cap C_+
%\end{array}
$$
Lerman proved that there are unique symplectic forms on
these orbifolds so
that $C_0$ is the symplectic reduction of $M$ at $0$,
the inclusions
$C_0\subset C_\pm$ are symplectic and the diffeomorphisms
$M_\pm - M_0 \to C_\pm - C_0$ are symplectomorphisms.
The symplectic orbifolds $C_+$ and $C_-$
are called the {\it symplectic cuts} of $M$ at $0$.

In this paper, we compute the rational
cohomology rings $H^*(C_+;\bbq)$ and $H^*(C_0;\bbq)$
in terms of those of
$M$ and $C_-$. In what follows,
the notation $H^*({\bf -})$ always means
$H^*({\bf -};\bbq)$.

We shall consider the kernel of the group homomorphism
$$H^*(C_-)\oplus H^*(M)\, \hfl{p^*-i^*}{} H^*(M_-)$$
(we denote by $p:M_\pm\to C_\pm$ the various projections and
by $i$ any inclusion between spaces). Note that $\ker (p^*-i^*)$ is a
subring of
$H^*(C_-)\oplus H^*(M)$. Consider the composed homomorphsim:
$$\begin{array}{lclc}
\delta : H^*(M_-,M_0) & \hfl{\hbox{\small diag}} &
H^*(M_-,M_0)\oplus H^*(M_-,M_0) & \hfl{\iso}\\
& \hfl{\iso} &
H^*(C_-,C_0)\oplus H^*(M,M_+) & \hfl{} \\
& \hfl{} & H^*(C_-)\oplus H^*(M). \end{array}$$
Observe that the image of $\delta$ is an ideal in $\ker (p^*-i^*)$.
Our main result is
\medskip\noindent\begin{ccote} {\bf Main Theorem\, }{\sl
The cohomology ring $H^*(C_+)= H^*(C_+;\bbq)$ is isomorphic, as a graded
ring, to
$\ker (p^*-i^*)/{Image\,}(\delta)$.}\end{ccote}

We now mention a few side results. First,
the cohomology ring of the symplectic reduction $C_0$ may be then obtained
from that of $C_-$ in the following way:

\begin{Proposition} \label{c-c0}
The cohomology ring $H^*(C_0)$ is the quotient of $H^*(C_-)$ by the
annihilator of the Poincar\'e dual of the suborbifold $C_0$. One has
a short exact sequence ($n=\dim M$):
$$0\to H_{n-*}(M_-;\bbq)\to H^*(C_-)\to H^*(C_0)\to 0.$$
\end{Proposition}

It is also worth noting that the cohomology ring of $M_-$ may be obtained
from those of $C_-$ and $C_0$. Let $i_! : H^{*-2}(C_0)\to H^*(C_-)$ be the
push-forward homomorphism.

\begin{Proposition} \label{shrick}
 There is  a short exact sequence
$$0\to H^{*-2}(C_0)\hfl{i_!}  H^*(C_-)\hfl{p^*} H^*(M_-)\to 0.$$
\end{Proposition}

As in \cite[Proposition 3.2]{HK} one can use
Propositions \ref{c-c0} and  \ref{shrick} to get the
Poincar\'e polynomials of
 $C_-$ and $C_0$ in terms of that of $M_-$:

\begin{Corollary} \label{betti}
The Poincar\'e polynomials for the rational cohomology of
 $C_-$ and $C_0$ are given by the equations
$$\begin{array}{lcl}
(1-t^2)P_{C_0}(t) & = & P_{M_-}(t) - t^n P_{M_-}(1/t) \\[2pt]
(1-t^2)P_{C_-}(t) &  = &
P_{M_-}(t) - t^{n+2} P_{M_-}(1/t).\end{array}$$
\end{Corollary}

\vskip .2 truecm\goodbreak
\noindent {\bf Remarks :}\nobreak %\vskip .2 truecm

\begin{ccote} \label{othcoeff}\rm
It seems likely that the results of this paper are true for the cohomology with
coefficients any ring in which the orders of the finite stabilizer groups
are invertible. The literature on cohomology of orbifolds
being still to be developed, we have followed,
in order to keep the paper focused on the problem at hand,
the lazy trend of using rational coefficients. (See also example
\ref{blowup}).
\end{ccote}

\begin{ccote} \label{equivar}  \rm
If $M$ admits a Hamiltonian action of a compact Lie group $G$ which
commutes with that of $S^1$ (equivalently: $f$ is $G$-invariant), the
the main theorem is also true for the equivariant cohomology
$H_G^*(C_+;\bbq)$. The proof is the same, using the $G$-equivariant
perfection of the Morse-Bott function $f$ (well known by experts). 
Note that the circle $S^1$ itself has this property.
\end{ccote}

\begin{ccote} \label{quantum} \rm
While we were pursuing this result, the paper \cite{IP}
appeared, which calculates some of the {\em quantum} cohomology of $M$
in terms of that of $M_+$, $M_-$, $C_0$, and $C_0$'s normal bundle in $C_+$.
Since the quantum cohomology subsumes the ordinary cohomology, it would be
interesting to see how our result is a consequence of those of \cite{IP}. 
\end{ccote}

\begin{ccote} \label{otherways} \rm
Since a symplectic cut can be defined as a symplectic reduction,
and there are now general techniques for computing the cohomology
of symplectic quotients (such as \cite{JK} for cohomology pairings
or \cite{TW} for generators and relations),
the question posed in this paper has other answers.
However, those answers are in terms of other ingredients
(like fixed point data in the other cut space, rather than its cohomology), 
so we believe that our result remains useful; 
in addition, the proofs are short and elementary 
(without e.g. heavy use of equivariant cohomology).
\end{ccote}

\section{Proofs}

{\bf Proof of the main theorem. }
Let $C$ denote the auxiliary space $C:= M/\sim \, = \, C_-\cup C_+$. 
Its cohomology will be calculated below. We will then make use of the maps
$M\dot\cup C_-\to C$ and $C_+\to C$ whose induced maps in cohomology 
will turn out to be respectively injective and surjective.

The long exact sequences of the pairs
$(C,C_-)$ and $(M,M_-)$ are linked by the commutative diagram
\begin{equation}\label{MVdiagr}\begin{array}{ccccccccccccccc}
\to & H^{*}(C,C_-) &\to & H^{*}(C)
&\hfl{i_-^*} & H^{*}(C_-) &\to & H^{*+1}(C,C_-) \to \\
& \downarrow\iso && \downarrow p^* &&\downarrow i_-^*
 &&  \downarrow\iso\\
\to & H^{*}(M,M_-) &\to & H^{*}(M)
&\hfl{i_-^*} & H^{*}(M_-) &\to & H^{*+1}(M,M_-) \to
%\\ &&&&&&
\end{array}\end{equation}
This produces a long exact sequence
\begin{equation}\label{longex1}
\to H^{*-1}(M_-) \, \hfl{\partial}  H^*(C)  \hfl{(i_-^*,p^*)}
H^*(C_-)\oplus H^*(M)
\hfl{p_-^*-i_-^*} H^*(M_-)  \to \end{equation}
where the connecting homomorphism $\partial$ is 
obtained using the vertical isomorphisms
$$\partial : H^{*-1}(M_-) \to H^{*}(M,M_-)
        \hfl{\iso} H^{*}(C,C_-) \to H^*(C).
$$
The moment map $f:M\to\bbr$ is a {\it perfect} Morse-Bott function
(See [Fr] when $M$ is K\"ahler; as observed by Atiyah
({\em Bull. London Math. Soc.} 14 (1982) p.7)
Frankel's proof works as well when $M$ symplectic).
This implies that $i^* :H^*(M)\to  H^*(M_-)$
is surjective. Therefore $\partial\equiv 0$ and
exact sequence \eqref{longex1} breaks:
\begin{equation}\label{longex1b}
0 \to  H^*(C)  \hfl{(i_-^*,p^*)}  H^*(C_-)\oplus H^*(M)
\hfl{p_-^*-i_+^*} H^*(M_-)  \to 0.\end{equation}
This first map is a ring homomorphism, induced from the map
$M\dot\cup C_- \to C$, locating $H^*(C)$ as a subring inside the known
ring $H^*(M)\oplus H^*(C_-)$.

Another consequence of the surjectivity of $i^* :H^*(M)\to  H^*(M_-)$
is, using diagram \eqref{MVdiagr}, that $i^* :H^*(C)\to  H^*(C_-)$
is also onto. We use the symmetric statement with $C_+$ and get the short exact
sequence:
\begin{equation} \label{exaseq44}
0\to H^*(C,C_+) \to H^*(C)\to  H^*(C_+) \to 0.
\end{equation}
Using the diagram
\begin{equation}\label{MVdiagr2}\begin{array}{ccccccccccccccc}
& 0 && 0 \\
& \downarrow && \downarrow \\
0\to & H^*(C,C_+) & \hfl{(i_-^*,p^*)} & H^*(C_-,C_0)\oplus H^*(M,M_+) \\
& \downarrow && \downarrow \\
0\to & H^*(C) & \hfl{(i_-^*,p^*)} & H^*(C_-)\oplus H^*(M)
\end{array}\end{equation}
and \eqref{exaseq44}  one sees that
$$H^*(C_+)\,\iso\, H^*(C)\big/H^*(C,C_+) \,\iso\,  \ker (p^*-
i^*)/{Image\,}(\delta) .$$
which proves our main theorem.
\cqfd

\vskip .5 truecm
{\bf Proof of Proposition \ref{c-c0}: }
It is already established in \eqref{exaseq44} that
 the cohomology exact sequence of the pair $C_-,C_0$ breaks:
$$0\to H^{*}(C_-,C_0)\to H^*(C_-)\to H^*(C_0)\to 0 .$$
 By excision and Poincar\'e
duality,
$$H^{*}(C_-,C_0) \iso H^{*}(M_-,M_0) \iso H_{n-*}(M_-).$$
The fact that $\ker \big( H^*(C_-)\to H^*(C_0)\big)$ is the
annihilator of the Poincar\'e dual of $C_0$ is proved as in
\cite[Proposition 3.3]{HK}. The argument there requires that
$i_* : H_*(C_0)\to H_*(C_-)$ is injective. This is guaranteed
by the perfection of the moment map $\bar f : C_-\to\bbr$ for the induced
$S^1$-action
on $C_-$.
\cqfd

\vskip .5 truecm
{\bf Proof of Proposition \ref{shrick}: }
Observe that we know from \eqref{exaseq44} that $p^*:H^*(C_-)\to H^*(M_-)$
is onto.
But Proposition \ref{shrick} actually comes from the following exact diagram:
$$\begin{array}{ccccccccccccc}
 & H^{*-2}(C_0) & \hfl{i_!} & H^*(C_-) & \hfl{p^*} & H^*(M_-) \\[2pt]
&&&&& {\scriptstyle\iso}\downarrow{\scriptstyle\cap [M_-]} \\[2pt]
& {\scriptstyle\iso}\downarrow
{\scriptstyle\cap [C_0]} && {\scriptstyle\iso}\downarrow
{\scriptstyle\cap [C_-]}&&
 H_{n-*}(M_-,M_0)\\
&&&&&\uparrow{\scriptstyle\iso}\\
0\to & H_{n-*}(C_0) & \hfl{i_*} & H_{n-*}(C_-) & \hfl{} &
H_{n-*}(C_-,C_0) &\to & 0.
\end{array}$$
The two zeros on the bottom line come from
the perfection of the moment map $\bar f : C_-\to\bbr$
for the induced $S^1$-action on $C_-$.
\cqfd

\section{Examples}

\begin{ccote}\label{blowup} \rm Suppose that the only critical point of $f$
in $M_-$
is an isolated minimum $u$.
Then $M_- = D^{2m}$, $C_0$ and $C_-$ are weighted projective spaces and $C_+$ is
 a ``weighted blowup" of $M$ (at $u$).

The relative cohomology $H^*(M,M_-;\bbz)$ is generated by the fundamental class
$v_- \in H^n(M_-,M_0)$ of $M_-$. Let $\delta(v)=(v_-,v)$.
We deduce from Corollary \ref{betti} that $C_-$ and $C_0$ have the
Poincar\'e polynomial of $\bbc P^{m}$ and $\bbc P^{m-1}$ respectively.
Since (to the top degree) all powers of their symplectic forms are nonzero,
the rings must be $H^*(M_-;\bbq )= \bbq [a]\big/(a^{m+1})$ and
$H^*(M_-;\bbq )= \bbq [a]\big/(a^{m})$. Passing from $\bbq$ to $\bbr$ for
the coefficients,
$a\in H^2(C_-;\bbr)$ may be chosen so that $a^m=v_-$. Our main theorem
then states that $H^*(C_+;\bbr)$ is the subring of
$$\bigg(\bbr [a]\big/(a^{m+1}) \oplus H^*(M;\bbr)\bigg) \bigg/ (a^m +v)$$
consisting of $\bbr (1,1)$ in degree 0, and everything in higher degree.
The above statement is presumably true with coefficients
$\bbz[w^{-1/m}]$ where $w$ is the product of the weights of the linear
action of $S^1$ on $T_uM$.
\end{ccote}

\begin{ccote} \rm Let $M$ be the space of Hermitian $3\times 3$ matrices
with eigenvalues
$-1,1,2$. This is a flag manifold, and acquires a symplectic form through
identification
with a coadjoint orbit of $U(3)$. Let $f$ be the (1,1)-matrix entry. Its
Hamiltonian flow
is periodic, being the conjugation by ${\rm diag}(e^{i\theta},1,1)$. The
critical
values of $f$ are again $-1,1,2$, each with a $\bbc P^1$ as critical set.
Let $A$ be the
minimum level set of $f$.

When we cut at $f=0$, we produce $M_-\sim A$ and $C_-\iso A \times \bbc P^2$.
The other cut $C_+$ is just a blowup of $M$ along $A$, and is in fact a
Bott-Samelson manifold \cite{BS}.

The cohomology ring $H^*(M)=H^*(M;\bbz)$ has as a $\bbz$-basis the Schubert
classes
(indexed by the Weyl group $\rm Sym_ 3$) \cite{Fu}:

\begin{tabbing}
\kern .5 truecm \= -- \ \= degree 0\ \=:\ \= $\sigma_{13}=1$ \\
\kern .5 truecm \= -- \ \= degree 2\ \=:\ \= $\sigma_{123}$, $\sigma_{132}$ \\
\>-- \> degree 4\ \>:\ \> $\sigma_{12}$, $\sigma_{23}$ \\
\>-- \> degree 6\ \>:\ \> $\sigma_{1}$
\end{tabbing}
Schubert calculus gives us the relators:

-- \ $\sigma_{123}\sigma_{132}=\sigma_{12}+\sigma_{23}$,

-- \ $\sigma_{123}^2=\sigma_{23}$ and $\sigma_{132}^2=\sigma_{12}$,

-- \ $\sigma_{123}\sigma_{12} = \sigma_{132}\sigma_{23}=\sigma_{1}$  and
$\sigma_{123}\sigma_{23} = \sigma_{132}\sigma_{12}=0$.

\vskip .3 truecm
Setting $u:=\sigma_{123}$ and $v:=\sigma_{132}$, one gets the presentation
$$H^*(M)\iso \bbz [u,v]\big/ (uv=u^2+v^2\, ,\, u^3=v^3=0\, , \,
u^2v=uv^2\, , \, u^2v^2=0)$$

Let $\bar a$ be the generator of $H^*(M_-) = \bbz[\bar a]/(\bar a^2=0)$. Set
also
$H^*(C_-) = \bbz[a,b]/(a^2=0,b^3=0)$. The homomorphisms
$p^*:H^*(C_-)\to H^*(M_-)$ and $i^*:H^*(M)\to H^*(M_-)$ are given by
$$p^*(a)=\bar a \, ,\, p^*(b)=0 \, ,\, i^*(u)=\bar a \, ,\, i^*(v)=0.$$
Setting $x:=(u,a)$, $y:=(v,0)$ and $z:=(0,b)$, the ring $H^*(C)$
is $\bbz [x,y,z]$ modulo

\ -- \ $xy=x^2+y^2$, $yz=0$

\ -- \ $x^2y=xy^2$ and  $x^3=y^3=z^3=0$.
\vskip .3 truecm
The image of $\delta : H^*(M_-,M_0) \to H^*(C_-)\oplus H^*(M)$
is generated by $(v^2,b^2)$ and $(uv^2,ab^2)$. Therefore, by our
main theorem,  the ring $H^*(M_+)$ is the quotient of $\bbz [x,y,z]$
by the following
relations:

\vskip .3 truecm
\ -- \ $xy=x^2+y^2$, $yz=0$ and $y^2=-z^2$

\ -- \ $x^2y=xy^2=-xz^2$ and  $x^3=y^3=z^3=0$.
\vskip .3 truecm

In this example, the cohomology ring $H^*(M_+)$ is known \cite[p. 993]{BS}:
$$H^*(M_+) = \bbz [a,b,c]\big/ \big(a^2,b(b+a), c(c-a+b)\big).$$
The isomorphism with our presentation is given by $a\mapsto y+z$,
$b\mapsto -z$ and $c\mapsto x+z$.
\end{ccote}

\vskip .3 truecm\goodbreak

\vskip .5 truecm\small
\noindent \parbox[t]{6 truecm}{Jean-Claude HAUSMANN\\
Math\'ematiques-Universit\'e\\ B.P. 240, \\
 CH-1211 Gen\`eve
 24, Suisse\\ hausmann@math.unige.ch} \ \hfill \hfill \
\parbox[t]{5 truecm}{Allen KNUTSON \\
 Department of Mathematics\\
 Brandeis University\\
 Waltham, MA 02254-9110 USA\\
 allenk@alumni.caltech.edu}

\end{document}